\documentclass[letterpaper,10pt,conference]{ieeeconf}

\usepackage{url}
\usepackage{tabularx}
\usepackage[font=small,labelfont=bf]{caption}
\usepackage[caption=false]{subfig}
\usepackage{scalerel}
\usepackage[utf8]{inputenc} 

\usepackage{amsmath,color,amsthm,amssymb,graphicx,amsfonts,float}
\usepackage{cite}

\newtheorem{prop}{Proposition}
\newcommand*{\Scale}[2][4]{\scalebox{#1}{$#2$}}%
\allowdisplaybreaks[2]

\usepackage{graphicx} 
\usepackage{amsfonts} 

\usepackage{booktabs} 
\usepackage{array} 
\usepackage{paralist} 
\usepackage{verbatim} 
\usepackage{subfig} 
\usepackage{multirow}
\usepackage{comment}
\usepackage{algorithm,algorithmic}
\usepackage[dvipsnames]{xcolor}
\usepackage{svg}
\IEEEoverridecommandlockouts
\newcommand{\yajing}[1]{{\color{blue}YL says: #1}}

\newcommand{\april}[1]{{\color{ForestGreen}AS says: #1}}

\newcommand{\andrey}[1]{{\color{RubineRed} AB says: #1}}

\title{Matrix Completion Using Alternating Minimization for Distribution System State Estimation}

\author{Yajing Liu, April Sagan, Andrey Bernstein, Rui Yang, Xinyang Zhou, and Yingchen Zhang
\thanks{Y. Liu, A. Bernstein, R. Yang, X. Zhou, and Y. Zhang are with Power System Engineering Center, National Renewable Energy Laboratory, Golden, CO 80401, {\small \{yajing.liu, andrey.bernstein, rui.yang, xinyang.zhou, yingchen.zhang\}@nrel.gov}}
\thanks{A. Sagan is with  the Department of  Mathematical Sciences, Rensselaer Polytechnic Institute, Troy, NY 12180, {\small sagana@rpi.edu}, and is supported by the National Science Foundation under Grant Number DMS-1736326.}
}

\begin{document}
\maketitle
\sloppy
\begin{abstract}
This paper examines the problem of state estimation in power distribution systems under low-observability conditions. The recently proposed constrained matrix completion method which combines the standard matrix completion method and power flow constraints has been shown to be effective in estimating voltage phasors under low-observability conditions using  single-snapshot information.  However, the method requires solving a semidefinite programming (SDP) problem, which becomes computationally infeasible for large systems and if multiple-snapshot (time-series) information is used. This paper proposes an efficient algorithm to solve the constrained matrix completion problem with time-series data. This algorithm is based on reformulating the matrix completion problem as a bilinear (non-convex) optimization problem, and applying the alternating minimization algorithm to solve this problem. This paper proves the summable convergence  of the proposed algorithm, and demonstrates its efficacy and scalability via IEEE 123-bus system and a real utility feeder system. This paper also explores the value of adding more data from the history in terms of computation time and estimation accuracy. 
\end{abstract}

\section{Introduction}
State estimation in power systems is critical for maintaining normal and secure operating conditions. In power transmission systems, the weighted least squares (WLS) method \cite{Abur2004} is  well developed and widely used by utilities for state estimation; however, the limited availability of measurements renders the WLS inapplicable in  power distribution systems. The requirement for state estimation in  distribution systems is becoming stringent because the increased penetration of distributed energy resources---such as solar photovoltaics, wind turbines, and energy storage systems---on distribution systems introduces bidirectional power flows that might impact the system responses to various types of disturbances.

Much work has been done in an attempt to address the low observability state estimation problem. Reference \cite{Baruki2017} determined the optimal minimum measurement locations and used the difference between the measured and calculated voltages and complex powers to obtain the voltage profile of the whole network. Reference \cite{Rigoni2017} used  the voltage sensitivity at each bus to the load at all buses to determine the voltage given the measured load at each location. References \cite{manitsas2012, Pertl2016,YCNN,zamzam2019data} applied neural networks, which do not require system models, to estimate the voltage in the distribution grid; however, these approaches still require the installation of significant numbers of phasor measurement units (PMUs) when considering large systems. References \cite{manitsas2012,Clements2011, wu2013,zhou2019gradient} used pseudo-measurements, which are estimates based on historical data to overcome the requirement of the installation of more meters; however, it is known that the pseudo-measurements typically have  larger measurement errors than the metered measurements.

An alternative approach that does not require installing new sensors or explicit computation of pseudo-measurements was recently proposed by \cite{Schmitt}. This method leverages the standard matrix completion method \cite{exactMatrixCompletion} and augments it with linearized power flow constraints to acknowledge for the physical network constraints.
This method was shown to be effective in estimating voltage phasors using single-snapshot information.  However, this method requires solving a semidefinite programming (SDP) problem, which becomes computationally infeasible for large systems or if multiple-snapshot (time-series) information is used.

In fact,  temporal correlation (dependency between measurements at different time steps) exists among data in addition to the spatial correlation (dependency between measurements at different locations) and the correlation among measurement types (characterized by  power flow equations). By modeling both the temporal and spatial correlation of the data, \cite{gao_wang2016, Liao2019} apply the standard matrix completion method to recover missing PMU data,  and \cite{Genes2019} applies the standard matrix completion method with Bayesian estimation to recover missing low-voltage distribution system data.


This paper proposes an efficient algorithm to solve the constrained matrix completion problem with time-series data. To this end, this paper formulates  matrix $X$ using time-series data  by considering both the  spatial  and  temporal correlation of the data. Similar to \cite{Schmitt}, this paper leverages a linear approximation of the power flow equations as constraints to characterize the correlation among different measurement types.
Then, this paper applies 
the alternating minimization method,  which was a winner in the Netflix Challenge \cite{Koren2009},  to solve the proposed matrix completion problem for state estimation in multiphase distribution systems.

To apply the alternating minimization method,  matrix $X$ is written in a bilinear form, i.e., $X=UV$, and then the algorithm finds the best $U$ and $V$ in an alternating fashion. Because of the bilinear term in the objective function, the problem becomes nonconvex with respect to any variable and thus the alternating minimization method is not guaranteed to converge to the global optima. However, 
results on the  convergence of the method were established under some regularity conditions in \cite{Luo1993,xu_yin_2013}. This paper proves that the objective function satisfies these conditions, which guarantee the convergence of the proposed method. 



The rest of this paper is organized as follows. Section~\ref{section:ProblemFormulate} introduces the standard matrix completion model, the formulation of the data matrix using time-series data, the linearized power flow constraints, and the full matrix completion model for state estimation in multiphase distribution systems. The alternating minimization algorithm and the  summable convergence are presented in Section~\ref{section:AM}. Simulation results on the IEEE 123-bus system and a real utility feeder are shown in Section~\ref{section:simulationresults}. Section~\ref{section:conclusion} concludes the paper.

\section{Problem Formulation}
\label{section:ProblemFormulate}
This section first reviews the standard matrix completion formulation, and then  introduces how to build the data matrix using the time-series data in multiphase distribution systems. Next, the linearized power flow constraints are introduced. Finally,  the full matrix completion model for state estimation in multiphase distribution systems is derived, which combines the standard matrix completion model with the linear power flow constraints.
\subsection {Standard Matrix Completion}
\label{subsection:SMC}
{Let $M\in\mathbb{R}^{m\times n}$  be the matrix we wish to reconstruct, and let $\Omega$ denote the set of indices for known elements of $M$. Define the observation operator ${P}_\Omega:\mathbb{R}^{m\times n} \rightarrow \mathbb{R}^{m\times n} $ as:
   $ {P}_\Omega(M)_{ij}=M_{ij}$ for $(i,j) \in \Omega$; 0 otherwise.

 The goal is to recover $M$ from ${P}_\Omega(M)$,  under the assumption that $M$ is low rank. Because  the rank function is nonconvex and the tightest convex  relaxation of  the  rank  function is the nuclear norm function  \cite{exactMatrixCompletion}, we consider the following regularized matrix completion model:
\begin{equation}
\begin{aligned}
\label{relaxedMC}
& \underset{X\in\mathbb{R}^{m\times n}}{\min}
& & ||X||_*+\frac{\mu}{2}||{P}_\Omega(X)-{P}_\Omega(M)||_F^2,
\end{aligned}
\end{equation}
where $X$ is the decision variable, $||X||_*=\sum_{i=1}^{\text{min}\{m,n\}} \sigma_i(X)$ is the nuclear norm of $X$ with $\sigma_i(X)$ denoting the $i${th} {largest} singular value of $X$, and $\mu>0$ is a weight parameter.

\subsection{Matrix Formulation  in Multiphase Distribution Systems}
\label{subsection:matrixformulation}
For the scope of this paper, 
we consider a general three-phase distribution network consisting of one slack bus and a given number of multiphase $PQ$ buses.

Matrix completion was shown  effective for state estimation in general three-phase distribution systems by \cite{Schmitt}, which formulated the data matrix $M$ using  single-snapshot information. Considering  the temporal correlation of the data, we set up the data matrix $M$ using time-series data. 

 Assume that the voltage phasor and other measurements at the slack bus are known. So we use only  measurements at nonslack buses to form the data matrix. {Let $\mathcal{P}$ denote the set of phases at all nonslack buses and $|\mathcal{P}|$ be the number of all phases at  nonslack buses.} 
 The measurements we  use in this matrix are the real and imaginary parts of voltage phasor, voltage magnitude, and the active and reactive power injection at each phase of nonslack buses.  Consider a time series $t=1,\ldots, T$. Let $M^t$ denote the measurement matrix at time $t$ {such that each column represents a phase and each row represents a quantity relevant to the phase. To be specific, for each phase $i\in\mathcal{P}$,
the corresponding column of $M^t$ is of the form:}
\[\begin{aligned}
\left[\Re(v_i),\ \Im(v_i),\ |v_i|,\ \Re(s_i), \ \Im(s_i)\right]^{\intercal},
\end{aligned}\]
where $\Re(\cdot)$ and $\Im(\cdot)$ are the real part and imaginary part of a complex variable, respectively, $\intercal$ is the transpose notation, $v=[v_1,\ldots, v_{|\mathcal{P}|}]^{\intercal}\in\mathbb{C}^{|\mathcal{P}|}$ is the vector containing voltage phasors at each phase of nonslack buses, and $s=[s_1,\ldots, s_{|\mathcal{P}|}]^{\intercal}\in\mathbb{C}^{|\mathcal{P}|}$ is the vector of power injections at each phase of nonslack buses. 
{Then the matrix $M$ is constructed by 
\begin{equation}
\label{eqn:dataMatrix}
\begin{aligned}
M=\left[M^1;\ M^2;\ \cdots; \ M^{T}\right] \in \mathbb{R}^{m \times n},
\end{aligned}
\end{equation}
where  $T$ is the length of the time  series, and $m = 5T, n=|\mathcal{P}|$.}
Note that the data matrix is not limited to the structure proposed above; it can accommodate any available measurements if these measurements are correlated such that $M$ has the (approximate) low-rank property.

Because the rows and/or columns of the matrix $M$ are related by temporal correlation, spatial correlation, or power flow equations, such a matrix has low rank. Fig.~\ref{singularvalues} shows that the low-rank property holds for the IEEE 123-bus  system with one-time step, two-time step, and three-time step  data formulation.

\begin{figure}
   \centering
   \includegraphics[scale = 0.42, trim= 46 0 0 0, clip]{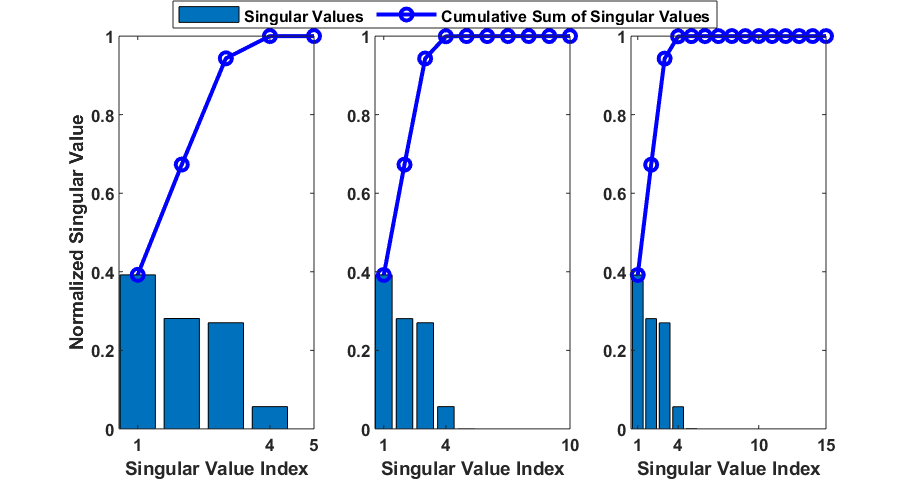}
   \caption{Singular  values  of  data  matrices  for  the  IEEE   123-bus  system with one-time step, two-time step, and three-time step  data matrix formulation. The data matrices are of size $5\times260$, $10\times260$, and $15\times260$, respectively. The  bars  show  the  individual  singular values (normalized by the singular value sum), and the circles show the cumulative sums of (normalized) singular values. 
   The first four largest singular values comprise 99\% of the singular value sum.} 
\label{singularvalues}
\end{figure}
  
\subsection{Linear Power Flow Constraints} \label{linearLoadFlow}

As in \cite{Schmitt},
we use the linear approximation  of  power  flows  as  constraints  to  characterize the correlation among different measurement types, which provides  more  information  about  the  system,  and  thus  likely  requires fewer  measurements  to  recover  the  matrix  than the standard  matrix  completion  formulation.


Recall that we use  $\mathcal{P}$, $v\in\mathbb{C}^{|\mathcal{P}|}$, and $s\in\mathbb{C}^{|\mathcal{P}|}$  to denote the set of phases, the vector of voltage phasors, and the vector of power injections at  all nonslack buses, 
respectively. 
 We employ approximations of the form:
\begin{subequations} \label{eq:lin_app}
\begin{align}
v  \; &\approx B \begin{bmatrix} \Re(s) \\ \Im(s) \end{bmatrix} + \: w,  \label{eq:vlin} \\
|v| &\approx C \begin{bmatrix} \Re(s) \\ \Im(s) \end{bmatrix}  +|w|, \label{eq:vmaglin}
\end{align}
\end{subequations}
where the coefficients $B,C\in\mathbb{C}^{|\mathcal{P}|\times 2|\mathcal{P}|}, w\in\mathbb{C}^{|\mathcal{P}|}$ are derived by the power flow linearization methods developed in, e.g., \cite{linearLoadFlow,Christakou2013,Giannakis2016}.
Write (\ref{eq:lin_app}) in the following equivalent form with $A_1,A_2,A_3,A_4,C_1,C_2\in\mathbb{R}^{|\mathcal{P}|\times |\mathcal{P}|}$:
\begin{subequations} \label{eq:lin_app2}
\begin{align}
\Re(v)  \; &\approx [A_1\quad A_2] \begin{bmatrix} \Re(s) \\ \Im(s) \end{bmatrix} + \: \Re(w),  \label{eq:vreallin} \\
\Im(v)  \; &\approx [A_3\quad A_4] \begin{bmatrix} \Re(s) \\ \Im(s) \end{bmatrix} + \: \Im(w),  \label{eq:vimglin} \\
|v| &\approx [C_1\quad C_2] \begin{bmatrix} \Re(s) \\ \Im(s) \end{bmatrix}  +|w|, \label{eq:vmaglin2}
\end{align}
\end{subequations}
where $\Scale[0.9]{[A_1 \quad A_2] = \Re(B),[A_3 \quad A_4]=\Im(B)}$, and $\Scale[0.9]{[C_1\quad C_2]=C}$.

 We use $v^t$, $s^t$ to denote the corresponding $v,  s$ in (\ref{eq:lin_app2}) at time slot $t$. Assume that $v_0$ (the slack bus voltage phasor) and the system topology are the same at different time steps, and then by using the linearization method in  \cite{linearLoadFlow}, $A_1, A_2, A_3, A_4, C_1, C_2$ are the same at different time steps.  Hence, we have the linear approximation of voltage phasor and magnitude at time $t$ in (\ref{eq:lin_app2}) with $v,s$ being replaced by $v^t$ and $s^t$.

 For simplicity of expression in the sequel, we use the following model to express this linear model at time $t=1,\ldots, T$: 
 \begin{equation}
\label{eqn:linmodel}
y\approx Ax+b,
\end{equation}
where: \[y = \begin{bmatrix}
    \Re(v^1)\\
    \Im(v^1)\\
    |v^1|\\
    \vdots\\
    \Re(v^{T})\\
    \Im(v^{T})\\
    |v^{T}|
\end{bmatrix},\  A = \begin{bmatrix}
    A_1 & A_2 & 0 & \cdots & 0 & 0\\
     A_3& A_4 & 0 & \cdots & 0 & 0 \\
   C_1& C_2 & 0 & \cdots & 0  &0\\
    \vdots & \vdots & \vdots & \ddots & 0 & 0 \\
    0 & 0 & 0 & \cdots & A_1 & A_2\\
    0 & 0 & 0 & \cdots & A_3 & A_4\\
    0 & 0 & 0 & \cdots & C_1 & C_2
\end{bmatrix},   \]
\[x= \begin{bmatrix}
    \Re(s^1)^{\intercal}\  \Im(s^1)^{\intercal}\  \ldots\  \Re(s^{T})^{\intercal}\ \Im(s^{T})^{\intercal}\end{bmatrix}^{\intercal},\] and
     \[b =\begin{bmatrix}
    \Re(w)^{\intercal}\  \Im(w)^{\intercal}\ |w|^{\intercal}\ \ldots\  \Re(w)^{\intercal}\ \Im(w)^{\intercal}\ |w|^{\intercal}
\end{bmatrix}^{\intercal}\]
with $ A\in\mathbb{R}^{3T|\mathcal{P}|\times 2T|\mathcal{P}|},x\in\mathbb{R}^{2T|\mathcal{P}|}$, and $y, b\in\mathbb{R}^{3T|\mathcal{P}|}$.

{The linear model (\ref{eqn:linmodel}) will be incorporated into the matrix completion formulation in the following subsection.}

\subsection{The Regularized Matrix Completion for State Estimation}

{
By incorporating the  linear power flow model (\ref{eqn:linmodel}) into  problem (\ref{relaxedMC}) as a regularization term, we obtain the following optimization problem:
 \begin{eqnarray}
 \hspace{-3mm}&{\min} & ||X||_*+\frac{\mu}{2}||{P}_\Omega(X-M)||_F^2 +\frac{\nu}{2}||y-(Ax+b)||_2^2 \nonumber\\
\hspace{-3mm}&\text{over} & X \in\mathbb{R}^{m\times n},y \in\mathbb{R}^{\frac{3}{5}mn},x \in\mathbb{R}^{\frac{2}{5}mn},\nonumber\\
\hspace{-3mm}&\text{s.t.} & \Scale[0.95]{y = \left[
    a_1^{\intercal}X \ \, 
    a_2^{\intercal}X\ \, 
    a_3^{\intercal}X\,
    \ldots \,
    a_{_{3T-2}}^{\intercal}X\ \,
    a_{_{3T-1}}^{\intercal}X\ \,
    a_{_{3T}}^{\intercal}X 
\right]^{\intercal}}, \nonumber\hspace{-3mm}\\
\hspace{-3mm}&&  \Scale[0.95]{x = \left[
   c_1^{\intercal}X\ \,
     c_2^{\intercal}X \,
    \ldots \,
    c_{_{2T-1}}^{\intercal}X\ \,
    c_{_{2T}}^{\intercal}X
         \right]^{\intercal}},\hspace{-3mm} \label{regularizedmc}
 \end{eqnarray} 
where, $\nu>0$ is a weight parameter,  $m = 5T, n = |\mathcal{P}|$,
$a_{3(t-1)+i}=e_{5(t-1)+i}$ and $c_{2(t-1)+i}=e_{5(t-1)+3+i}$ are standard basis vectors in $\mathbb{R}^{m}$, in which $t$ denotes the time step number and $i$ denotes the $i$th $a$ or $c$ at the $t$-th time step. For example, if $a_1$, $t =1, i =1$, then, $a_1 = e_1\in\mathbb{R}^m$; if $c_1$, $t=1, i = 1$, then $c_1=e_4\in\mathbb{R}^m$. 

\section{Alternating Minimization for Solving the Full Matrix Completion Model}
\label{section:AM}
This section first  proposes the alternating minimization algorithm for the full matrix completion model in multiphase distribution system state estimation, and then proves the  summable convergence of the proposed algorithm.
\subsection{Alternating Minimization Algorithm}
Any matrix $X\in\mathbb{R}^{m\times n}$ of a rank up to $r$ has a matrix product $X= UV$ form with $U\in\mathbb{R}^{m\times r}$ and $V\in\mathbb{R}^{r\times n}$, and its nuclear norm can be expressed by the Frobenius norm of $U$ and $V$ as follows \cite{Srebro2004}:
{
\begin{equation}
\label{eqn:nucleareqFrobenius}
\begin{aligned}
||X||_*:=\ &\underset{U \in \mathbb{R}^{m \times r},V \in \mathbb{R}^{r \times n}}{\arg\min}\ \frac{1}{2} (||U||_F^2+||V||_F^2) \\
& \quad\ \ \mathrm{s.\,t.~} \quad\quad\quad\ 
 X=UV.
\end{aligned}
\end{equation}
}%
Note that for sufficiently small $r$ ($r\ll \text{min}\{m,n\}$), using this characterization of the nuclear norm  allows us to dramatically reduce the size of the problem.  We now use (\ref{eqn:nucleareqFrobenius}) to reformulate (\ref{regularizedmc}) as follows:
\begin{eqnarray}
\hspace{-4mm}& \Scale[0.95]{\underset{\substack{U \in \mathbb{R}^{m \times r},\\V \in \mathbb{R}^{r \times n}}}{\min}}&\hspace{-2mm}
\Scale[0.95]{\frac{1}{2} (||U||_F^2+||V||_F^2)+\frac{\mu}{2} ||{P}_\Omega(UV-M)||_F^2}\nonumber\\[-12pt]
\hspace{-4mm}&&\hspace{-2mm}\quad\quad \quad\quad \Scale[0.95]{+\frac{\nu}{2}||f_1(UV)-(Af_2(UV)+b)||_2^2,}\nonumber\\
\hspace{-4mm}&\mathrm{s.t.}&\hspace{-4mm} \Scale[0.86]{f_1(UV) = \left[
    a_1^{\intercal}X \ \, 
    a_2^{\intercal}X\ \, 
    a_3^{\intercal}X \,
    \ldots \,
    a_{_{3T-2}}^{\intercal}X\ \,
    a_{_{3T-1}}^{\intercal}X\ \,
    a_{_{3T}}^{\intercal}X 
\right]^{\intercal}}, \nonumber\\
\hspace{-4mm}&&\hspace{-4mm} \Scale[0.86]{f_2(UV) = \left[
   c_1^{\intercal}X\ \,
     c_2^{\intercal}X \,
    \ldots \,
    c_{_{2T-1}}^{\intercal}X\ \,
    c_{_{2T}}^{\intercal}X
         \right]^{\intercal}},\nonumber\\
\hspace{-4mm}&&\hspace{-4mm} \Scale[0.86]{X =UV.}\label{bilinearmc}
\end{eqnarray}

The alternating minimization algorithm that updates the variables in an alternating fashion is proved to be one of the most accurate and effective methods for solving the matrix completion problem \cite{altminmatcomp}. For our matrix completion formulation in state estimation, the alternating minimization algorithm solves for $U$ and $V$ alternatively while fixing the other factor.
The pseudo-code of the alternating minimization algorithm for the formulation (\ref{bilinearmc}) is given in Algorithm~\ref{Agorthim:AM}.

\begin{algorithm}
 \caption{Alternating Minimization Algorithm for Matrix Completion in Distribution System State Estimation}
 \begin{algorithmic}[1]
 \renewcommand{\algorithmicrequire}{\textbf{Input:}}
 \renewcommand{\algorithmicensure}{\textbf{Initialize:}}
 \REQUIRE Set of known indices $\Omega$, measurement matrix $M$, linear model ${A}, b$, and the number of iterations $N$. 
\ENSURE Calculate the rank $r$ singular value decomposition  $[\tilde{U} \tilde{\Sigma} \tilde{V}]=\text{SVD}_r({P}_\Omega(M))$, and set $V^{(0)}={{\tilde{\Sigma}}^{\frac{1}{2}}}\tilde{V}$.
  \FOR {$k = 1,\ldots,N$ 
  }
\STATE \[
\begin{aligned}
\label{eqn_U_update_centralized}
U^{(k)} &= \underset{\substack{U \in \mathbb{R}^{m \times r}}}{\arg\min}\
||U||_F^2 +\mu ||{P}_\Omega(U{V^{(k-1)}}-M)||_F^2\\
& +{\nu}||f_1(U{V^{(k-1)}})-(Af_2(U{V^{(k-1)}})+b)||_2^2
\end{aligned}
\]

\STATE \[
\begin{aligned}
\label{eqn_V_update_centralized}
V^{(k)}&=\underset{\substack{V \in \mathbb{R}^{r \times n}}}{\arg\min}\
||V||_F^2+{\mu} ||{P}_\Omega(U^{(k)}{V}-M)||_F^2\\
&+{\nu}||f_1(U^{(k)}V)-(Af_2(U^{(k)}{V})+b)||_2^2
\end{aligned}\]
  \ENDFOR
 \RETURN $X=U^{(N)}{V^{(N)}}$ 
 \end{algorithmic} 
 \label{Agorthim:AM}
 \end{algorithm}
\subsection{Convergence of the Alternating Minimization Algorithm}


Because of the bilinear form $X=UV$, problem~(\ref{bilinearmc}) is not jointly convex with respect to any variable.
In general, there is no guarantee that a stationary point of a nonconvex problem will coincide with the global optimum. However, \cite{xu_yin_2013} proved the  summable convergence of the alternating minimization algorithm if the objective function satisfies the following four conditions: 1) it  is  continuous  in  its  domain  and  always bigger  than  negative  infinity;  2) it  has  a Nash equilibrium; 3) it is strongly convex with respect  to  each  variable; and 4) it satisfies the Kurdyka-Łojasiewicz property (please see Definition~2.5 in  \cite{xu_yin_2013} for the definition). We show below that our objective function satisfies these four conditions;  thus, we have the following proposition. 

\begin{prop}
The alternating minimization algorithm (Algorithm~1) satisfies the summable convergence, i.e.: \[\sum_{k=1}^{+\infty}||U^{(k+1)}-U^{(k)}||_F+||V^{(k+1)}-V^{(k)}||_F<+\infty.\] 
\end{prop}

\begin{proof}
The objective function $f$ can be written as follows:
\begin{align}
    \label{objfun}
    f(U,V) = &\frac{1}{2} (||U||_F^2+||V||_F^2)+\frac{\mu}{2} ||{P}_\Omega(UV-M)||_F^2\nonumber\\
&\ +\frac{\nu}{2}||f_1(UV)-(Af_2(UV)+b)||_2^2,
\end{align}
where:
\begin{align}
    \label{lastterm}
    &||f_1(UV)-(Af_2(UV)+b)||_2^2\nonumber\\
    =&\Scale[0.9]{\sum\limits_{i=1}^{T}||X^{\intercal}a_{_{3i-2}}-(A_1X^{\intercal}c_{_{2i-1}}+A_2X^{\intercal}c_{_{2i}})-\Re(w)||_2^2}\nonumber\\
    &\ \ \Scale[0.9]{+\sum\limits_{i=1}^{T}||X^{\intercal}a_{_{3i-1}}-(A_3X^{\intercal}c_{_{2i-1}}+A_4X^{\intercal}c_{_{2i}})-\Im(w)||_2^2}\nonumber\\
    &\ \ \Scale[0.9]{+\sum\limits_{i=1}^{T}||X^{\intercal}a_{_{3i}}-(C_1X^{\intercal}c_{_{2i-1}}+C_2X^{\intercal}c_{_{2i}})-|w|||_2^2}
\end{align}
with $X=UV$.

First, the domain of the objective function is $(\mathcal{U},\mathcal{V})=\{(U,V):U\in\mathbb{R}^{m\times r}, V\in\mathbb{R}^{n\times r}\}$.
By the form of $f(U,V)$ in (\ref{objfun}), we have that $f(U,V)$ is continuous in its domain $(\mathcal{U},\mathcal{V})$, and it is non-negative for any $(U,V)\in (\mathcal{U},\mathcal{V})$ .

Second,  to prove that (\ref{bilinearmc}) has a Nash point, we prove that (\ref{bilinearmc}) has a global minimizer, which by definition is a Nash point.  By  the definition of $f(U,V)$ in (\ref{objfun}), the following expression holds:
\[\lim\limits_{\substack{||U||\rightarrow +\infty,\\
||V||\rightarrow +\infty}} f(U,V)=+\infty,\]
which means that $f(U,V)$ is coercive \cite{Lamberscoercive}.
 Then, by Theorem~2 in \cite{Lamberscoercive}, problem (\ref{bilinearmc}) has a global minimizer.

Third, we prove that $f(U,V)$ is strongly convex with respect to $U$ and $V$, respectively. 
To prove that $f(U,V)$ is strongly convex with respect to $U$, we prove that $\frac{\partial^2f(U,V)}{\partial U^2}\succ 0$ as follows.

Let $g(U):=||f_1(UV)-(Af_2(UV)+b)||_2^2$ for a fixed $V$. Then we have 

\begin{align}
    \label{secondorderder}
    \frac{\partial^2f(U,V)}{\partial U^2}= &I_{mr\times mr}+\mu (V\otimes I_{m\times m})^T {R_\Omega} (V\otimes I_{m\times m}) \nonumber\\ &\quad\quad+\frac{\nu}{2}\frac{\partial^2g}{\partial U^2},
\end{align}
where $R_\Omega\in\mathbb{R}^{mr\times mr}$ is the diagonal matrix such that 
$(R_\Omega)_{ii}=1$ if  $\big( ( i \mod n) +1, \lceil \frac{i}{m} \rceil \big) \in \Omega$ and 0 otherwise.

By the form of $g$ in (\ref{lastterm}), we  focus only at the second-order partial derivative of the first term in the first summation term in (\ref{lastterm}); treatment of the other terms can be derived similarly. Let $g_1$  denote the first term in the first summation term in (\ref{lastterm}).
We will prove that 
$g_1$ is a convex function by definition, then using that $g_1:\mathbb{R}^{m\times r}\rightarrow \mathbb{R}$ is convex if and only if $\frac{\partial^2g_1}{\partial U^2}\succeq 0$ to derive that  $\frac{\partial^2g_1}{\partial U^2}\succeq 0$.

By (\ref{lastterm}), we have that:
\[g_1(U) = ||A_1V^{\intercal}U^{\intercal}c_1+A_2V^{\intercal}U^{\intercal}c_{2}-V^{\intercal}U^{\intercal}a_1+\Re(w)||_2^2.\]
Because $\mathcal{U}$ denotes the set of all matrices of size $m\times r$, we have that $\mathcal{U}$ is a convex set.
Then by the form of $g_1(U)$, it is easy to derive that  $\forall U_1, U_2\in \mathcal{U}, \forall t\in [0,1]$: \[g_1(tU_1+(1-t)U_2)\leq tg_1(U_1)+(1-t)g_1(U_2).\]

Hence, $g_1$ is convex, implying that $\frac{\partial^2g_1}{\partial U^2}\succeq 0$.

Because all the other terms in (\ref{lastterm}) share the similar form with $g_1$, we can similarly prove that all the other terms in (\ref{lastterm}) are convex functions with respect to $U$; thus, we have that $\frac{\partial^2g}{\partial U^2}\succeq 0.$
By (\ref{secondorderder}), we have that 
$I_{mr\times mr}\succ 0$, and $\mu (V\otimes I_{m\times m})^T {R_\Omega} (V\otimes I_{m\times m}) \succeq 0$; therefore, for $\mu>0, \nu>0, \frac{\partial^2g}{\partial U^2}\succeq 0$, we have
$\frac{\partial^2f(U,V)}{\partial U^2}\succ 0.$

Hence, we conclude that $f$ is strongly convex with respect to $U$. The proof that $f$ is strongly convex with respect to $V$ is similar and is omitted for brevity.

By \cite{xu_yin_2013}, any real analytic function satisfies the Kurdyka-Łojasiewicz  property. And by definition of the real analytic function, the objective function (\ref{objfun}) is real analytic; thus, (\ref{objfun}) satisfies the Kurdyka-Łojasiewicz property.

This concludes the proof.
\end{proof}
\section{Simulation Results}
\label{section:simulationresults}
This section demonstrates Algorithm~\ref{Agorthim:AM} on the IEEE 123-bus test system \cite{ieee123} and  a real utility feeder system. The two systems are both three-phase unbalanced radial distribution systems, in which buses are one-, two-, or three-phase. The slack bus is three-phase, and the total number of phases at all buses for the IEEE 123-bus system is 263 and that for the real utility feeder system is 1234. The data for both systems were simulated at 1-minute resolution using power flow analysis with diversified load and solar profiles that were created for each bus using realistic solar irradiance and load consumption data.
The voltage magnitudes range from 0.95 to 1 p.u., and the voltage angles are around 0 or $\pm 120$ degrees.
\vfill
\begin{figure}
   \centering
   \includegraphics[scale=0.45,trim= 80 240 45 240, clip]{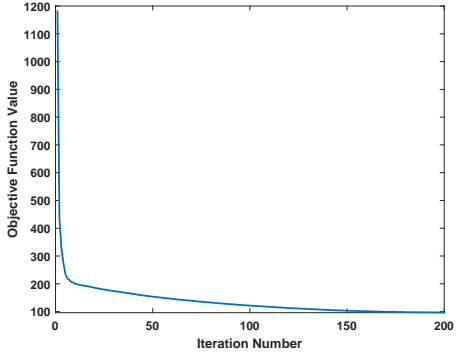}
   \caption{Convergence of the alternating minimization algorithm with one-time step data matrix formulation and 50\% available measurements.}
\label{convergence}
\end{figure}
\vfill
\begin{figure}
   \centering
   \includegraphics[scale=0.5,trim= 95 210 50 200, clip]{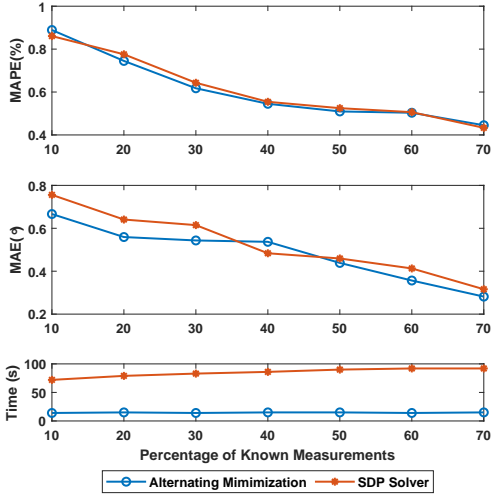}
   \caption{Comparison of performance for the alternating minimization algorithm and SDP solver. }
\label{mape_mae_onetime}
\end{figure}

As in Section~\ref{subsection:matrixformulation}, the data matrix is formulated using  the  real and imaginary parts of voltage phasor,  voltage  magnitude,  and active  and reactive power injection at each phase of nonslack buses. We consider  the real and imaginary parts of voltage phasor as variables, and the voltage magnitude,  active power, and reactive power as potentially known measurements. In our experiments, we define the percentage of known measurements as the ratio of the number of known measurements and the total number of voltage magnitude, active power, and reactive power in the data matrix. If less than two thirds of the measurements are known, there are fewer measurements than variables, under which case, the system is underdetermined and the WLS cannot generate a unique solution. For all of our experiments,  the known measurements are added with noise with zero mean and 1\% of the true value as the  standard deviation. And the simulation results are based on the average of $50$ runs for each scenario.
\subsection{Performance on IEEE 123-Bus System}
First, we implement our algorithm with one-time step data formulation at one time slot when solar generations are nonzero with $10\%$ to $70\%$ measurements available (where WLS is under-determined). We compare
the alternating minimization algorithm with the SDP solver used by \cite{Schmitt}. Both algorithms are run with MATLAB on a laptop with 1.9-GHz CPU and 32 GB RAM. 

 Fig.~\ref{convergence} shows one representative convergence result (with 50\% available measurements)  of the alternating minimization algorithm. 
Fig.~\ref{mape_mae_onetime} shows that the mean absolute percentage errors (MAPEs) of the voltage magnitude  and the mean absolute errors (MAEs) of the voltage angle are comparable for the  alternating minimization algorithm and the SDP solver. The running time for the  alternating minimization algorithm is about 15 seconds  for different percentages of available measurements, whereas that for the SDP solver is more than 70 seconds. These results show that the  alternating minimization algorithm estimates voltage phasors more  computationally efficient than the SDP solver under similar performance. The running time for the SDP solver increases  as more measurements are available. This is  because the SDP solver solves a constrained SDP problem, where more measurements being available means more constraints to be satisfied, and thus longer time is required to achieve the convergence.

\vfill
\begin{figure}
   \centering
   \includegraphics[scale=0.5, trim= 95 200 50 190, clip]{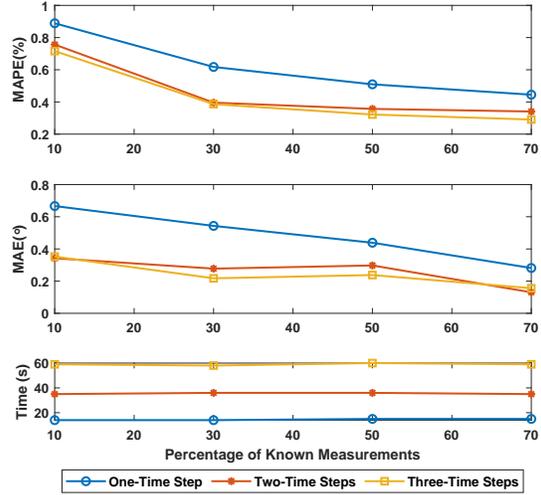}
   \caption{Comparison of performance for alternating minimization algorithm on one-time step, two-time steps, and three-time steps. }
\label{mape_mae_timeseries}
\end{figure}
To test the computation  time  and  estimation  accuracy as more data are used, we  implement the  alternating minimization algorithm on the two-time step and three-time step data matrix formulation with $10\%, 30\%, 50\%$, and  $70\%$ measurements available. Fig.~\ref{mape_mae_timeseries} shows that the MAPEs of the voltage magnitudes and the MAEs of the voltage angles for the two-time step and three-time step cases are much better than the one-time step case. This is because the matrix completion problem minimizes the rank, and the larger the time steps are, the smaller the rank of the matrix is compared to its size (which is shown in Fig.~\ref{singularvalues}) and the more accurate the estimates are. The  running time also increases with the length of time steps, however, because of the increased size of the matrices.
The MAPEs of the voltage magnitudes and MAEs of the voltage angles for 
the three-time step case is slightly better than the two-time step case overall, but the running time almost doubles. However, the running time for the three scenarios shows that the alternating minimization algorithm for state estimation can be implemented in real time for the IEEE 123-bus system.


\subsection{Performance on A Real Feeder System}
In order to demonstrate the scalability of our algorithm, we implement it on a real utility feeder system with 1234 phases under the assumption that 50\%  measurements are available  for one-time, two-time, and three-time step matrix formulation, respectively. The performance in terms of MAPE, MAE, and running time are shown in Table~\ref{tab:realfeeder}. The results show that the performance improves with more data being used for building the matrix and the computation time also increases. However, the running time shows that the algorithm can be applicable in real-time estimation.
\begin{table}
	\caption{Performance  on a real utility feeder system}\vspace{-10pt}
\label{tab:realfeeder}
	\begin{center}
	\begin{tabular}{l c c c c} 
			\toprule
	            & MAPE  & MAE &  Time (s)\\
			\midrule
			One-time step matrix formulation &  $0.4840$ & $0.3967$ & $50$\\
			Two-time step matrix formulation &  $0.4240$ & $0.2869$ & $105$\\
			Three-time step matrix formulation & $0.3503$ & $0.2265$ & $156$ \\
			\bottomrule
		\end{tabular}\vspace{-15pt}
	\end{center}
\end{table}
\section{Conclusions}
\label{section:conclusion}
This paper applied the constrained matrix completion formulation for state estimation in multiphase distribution systems. Different from the constrained matrix completion, the data matrix was formulated using time-series data and  the alternating minimization algorithm was applied to solve the  constrained matrix completion. The  summable convergence of the alternating minimization algorithm was proved. In addition, the efficacy and scalability of the algorithm were demonstrated via the IEEE 123-bus system and a real utility feeder system.  

\bibliographystyle{IEEEtran}
\bibliography{references}

\end{document}